\documentclass[11pt]{article}
\usepackage{amsmath,amssymb,amsthm}

\newtheorem{theo}{Theorem}[section]
\newtheorem{defi}[theo]{Definition}

\newtheorem{lemm}[theo]{Lemma}
\newtheorem{corol}[theo]{Corollary}

\newcommand{\Out}{\operatorname{Out}}
\newcommand{\Ad}{\operatorname{Ad}}
\newcommand{\Char}{\operatorname{Char}}
\newcommand{\Aut}{\operatorname{Aut}}
\newcommand{\cF}{\mathcal{F}}
\newcommand{\cL}{\mathcal{L}}
\newcommand{\id}{\mathord{\text{\rm id}}}
\newcommand{\FAlg}{\operatorname{FAlg}}
\newcommand{\si}{\sigma}
\newcommand{\SL}{\operatorname{SL}}
\newcommand{\Z}{\mathbb{Z}}
\newcommand{\recht}{\rightarrow}
\newcommand{\actson}{\curvearrowright}
\newcommand{\N}{\mathbb{N}}
\newcommand{\cA}{\mathcal{A}}
\newcommand{\cG}{\mathcal{G}}
\newcommand{\QN}{\operatorname{QN}}
\newcommand{\B}{\operatorname{B}}
\newcommand{\Tr}{\operatorname{Tr}}
\newcommand{\ot}{\otimes}
\newcommand{\ox}{\overline{x}}
\newcommand{\Irred}{\operatorname{Irred}}
\newcommand{\R}{\mathbb{R}}
\newcommand{\M}{\operatorname{M}}
\newcommand{\C}{\mathbb{C}}
\newcommand{\Rep}{\operatorname{Rep}}
\newcommand{\grp}{\operatorname{grp}}
\newcommand{\cH}{\mathcal{H}}
\newcommand{\cU}{\mathcal{U}}
\newcommand{\al}{\alpha}
\newcommand{\om}{\omega}
\newcommand{\cGr}{\cG^\circ}
\newcommand{\Gammar}{\Gamma^\circ}
\newcommand{\cI}{\mathcal{I}}
\newcommand{\Q}{\mathbb{Q}}

\begin{document}
\begin{center}
{\LARGE\bf Every compact group arises as the outer \vspace{0.5ex}\\ automorphism group of a II$_1$ factor}

\bigskip

{\sc by S\'{e}bastien Falgui\`{e}res and Stefaan Vaes\footnote{Both authors are partially supported by Research Programme G.0231.07 of the Research Foundation -- Flanders (FWO).}}
\end{center}

{Department of Mathematics, K.U.Leuven, Celestijnenlaan 200B, B-3001 Leuven, Belgium \\
E-mail: sebastien.falguieres@wis.kuleuven.be, stefaan.vaes@wis.kuleuven.be}

\begin{abstract}
\noindent We show that any compact group can be realized as the outer automorphism group of a factor of
type II$_1$. This has been proved in the abelian case by Ioana, Peterson and Popa \cite{IPP} applying Popa's
deformation/rigidity techniques to amalgamated free product von Neumann algebras. Our methods are a
generalization of theirs.
\end{abstract}

\section*{Introduction}

The outer automorphism group $\Out(M)$ of a II$_1$ factor $M$ provides, in principle, a useful invariant to distinguish between families of II$_1$ factors. But this group $\Out(M)$ is extremely hard to compute.

Breakthrough rigidity results in the theory of II$_1$ factors were obtained recently by Popa (see \cite{sm1,sm2,popa1}) and are based on Popa's deformation/rigidity technique. Techniques and results of Popa were applied in several papers (see also \cite{bourbaki} for an overview) and two papers included as an application complete computations of outer automorphism groups of certain II$_1$ factors.
\begin{itemize}
\item In \cite{IPP}, it is shown that there exists, for every \emph{compact abelian} group $K$, a type II$_1$ factor $M$ with $\Out(M) \cong K$. The result in \cite{IPP} is an existence theorem and answered in particular the long standing open problem on the existence of II$_1$ factors without outer automorphisms.
\item In \cite{Popa Vaes}, type II$_1$ factors $M$ with $\Out(M)$ any \emph{discrete group of finite presentation} are explicitly constructed. This in particular gave the first explicit examples of II$_1$ factors without outer automorphisms.
\end{itemize}
In our paper, the methods of \cite{IPP} are applied to prove the existence of II$_1$ factors $M$ such that $\Out(M)$ is any, possibly non-abelian, compact group. In fact, given a minimal action of a compact group $G$ on the hyperfinite II$_1$ factor $R$, we prove the existence of an action of $\Gamma = \SL(3,\Z)$ on the fixed point algebra $R^G$ such that the II$_1$ factor $M$ given as the amalgamated free product $M = (R^G \rtimes \Gamma) *_{R^G} R$ satisfies $\Out(M) \cong G$.

The first rigidity results for II$_1$ factors are due to Connes in \cite{connes}, where it is shown in particular that $\Out(N)$ is a countable group whenever $N=\cL(\Gamma)$ is the group von Neumann algebra of an ICC property (T) group $\Gamma$. So, for concrete ICC property (T) groups $\Gamma$, the group $\Out(\cL(\Gamma))$ is in principle computable, although we do not know of any explicit computation.

Type II$_1$ factors admit a more general type of symmetry, under the form of \emph{finite index bimodules}. The finite index $M$-$M$-bimodules (modulo isomorphism) form a \emph{fusion algebra} that we denote as $\FAlg(M)$. Such a fusion algebra is a set equipped with an additive (direct sum) and a multiplicative (tensor product) structure and in which every element is the finite direct sum of irreducible elements. Another generalization of \cite{IPP} was given by the second author in \cite{vaes}, where the existence of II$_1$ factors $M$ with trivial fusion algebra, was shown. In our paper the fusion algebra $\FAlg(R)$ of the hyperfinite II$_1$ factor plays an important role: we make use of the fact that two countable fusion subalgebras of $\FAlg(R)$ become free after conjugating one of them by a well chosen automorphism of $R$ (see \cite{vaes}).

\section{Preliminaries}

We denote by $(M,\tau)$ a von Neumann algebra $M$ equipped with a faithful normal tracial state $\tau$. We denote
$M^n := M_n(\C) \otimes M$ for all $n \in \N$.

Let $(M,\tau)$ be a tracial von Neumann algebra and $N \subset M$ a von Neumann subalgebra. The $^*$-algebra of elements quasi-normalizing $N$ is defined as
\[\QN_M(N):=\{a \in M \mid  \exists a_1,\dots, a_n,b_1,\dots,b_m \in M\ \text{such that}\ Na \subset
\sum_{i=1}^n a_i N\ \text{and}\ aN \subset \sum_{i=1}^m N b_i \}\] The inclusion $N \subset M$ is called
\emph{quasi-regular} if $\QN_M(N)'' = M$.

If $N$ is a von Neumann subalgebra of a von Neumann algebra $M$, we denote by $\Aut(N \subset M)$ the group of
automorphisms of $M$ leaving $N$ globally invariant.

\subsection{Amalgamated free products}

We make use of amalgamated free product factors and recall some basic facts and
notations (see \cite{popa-amal} and \cite{voi} for more details). Let $(M_0, \tau_0)$ and $(M_1,\tau_1)$ be tracial von Neumann algebras
with a common von Neumann subalgebra $N$ such that ${\tau_0}_{|N}={\tau_1}_{|N}$. We denote by $E_i$ the unique $\tau$-preserving conditional expectation of $M_i$ onto $N$. The amalgamated free product $M_0 *_N M_1$ is, up to $E$-preserving isomorphism, the unique pair $(M,E)$ satisfying the following two conditions.
\begin{itemize}
\item The von Neumann algebra $M$ is generated by embeddings of $M_0$ and $M_1$ that are identical on $N$, and is equipped with a conditional expectation $E : M \recht N$.
\item The subalgebras $M_0$ and $M_1$ are free with amalgamation over $N$ with respect to $E$. This means that $E(x_1 \cdots x_n) = 0$ whenever $x_j \in M_{i_j}$ such that $E_{i_j}(x_j)=0$ and $i_1 \neq i_2$, \linebreak $i_2 \neq i_3, \ldots, i_{n-1} \neq i_n$.
\end{itemize}

The amalgamated free product $M_0 *_N M_1$ has a dense $^*$-subalgebra given by:
$$N \;\;\oplus\;\; \bigoplus_{n \geq 1} \;\; \left(
\bigoplus_{i_1 \neq i_2, \ldots, i_{n-1} \neq i_n} \overset{\circ}{M_{i_1}}
\cdots \overset{\circ}{M_{i_n}} \right)$$
where $\overset{\circ}{M_{i_k}} := M_{i_k}\ominus N$. The von Neumann algebra
$M_0 *_N M_1$ has a trace, defined by $\tau:={\tau_0}\circ E = {\tau_1}\circ E $.

\subsection{Popa's intertwining-by-bimodules technique}

In this paper, we use Popa's intertwining-by-bimodules technique (see \cite{sm1,sm2,popa1}) that we briefly recall now. Let $(B,\tau)$ be a tracial von Neumann algebra and $H_B$ a right Hilbert $B$-module. There exists a projection $p \in \B(\ell^2(\N)) \otimes B$ such that $H_B \cong p (\ell^2(\N) \ot L^2(B,\tau))_B$ and this projection is uniquely defined up to equivalence of projections in $\B(\ell^2(\N)) \otimes B$. We denote $\dim(H_B) := (\Tr \ot \tau)(p)$. Observe that the number $\dim(H_B)$ depends on the choice of tracial state $\tau$ in the non-factorial case.

Suppose now that $A$ and $B$ are two possibly non-unital von Neumann subalgebras of a tracial von Neumann algebra
$(N,\tau)$. We say that $A$ embeds into $B$ inside $N$, and write $A \prec_N B$, if there exists a non-zero $A$-$B$-subbimodule $H$ of $1_A L^2(N) 1_B$ such that $\dim(H_B) < + \infty$. The relation $A \prec_N B$ is independent of the choice of tracial state $\tau$ and is equivalent with the existence of $n \in \N$, $v \in \M_{1,n}(\C) \ot 1_A N 1_B$ a non-zero partial isometry and $\psi : A \recht \M_n(\C) \ot B$ a possibly non-unital $^*$-homomorphism satisfying $a v = v \psi(a)$ for all $a \in A$. For details, we refer to Section 2 in \cite{sm1} and Appendix~C in \cite{bourbaki}.

\subsection{Fusion algebras}

We first recall the abstract notion of a \emph{fusion algebra} and give below the basic example of the fusion algebra $\FAlg(M)$ of finite index bimodules over the II$_1$ factor $M$.

\begin{defi}
A fusion algebra $\cA$ is a free $\N$-module $\N[\cG]$ equipped with the following additional structure:
\begin{itemize}
\item an associative and distributive product operation, and a multiplicative unit element $e \in \cG$,
\item an additive, anti-multiplicative, involutive map $x \mapsto \overline{x}$, called \emph{conjugation},
\end{itemize}
satisfying Frobenius reciprocity: defining the numbers $m(x,y;z) \in \N$ for $x,y,z \in \cG$ through the formula
$$x y = \sum_z m(x,y;z) z$$
one has $m(x,y;z) = m(\ox,z;y)$ for all $x,y,z \in \cG$.
\end{defi}
The base $\cG$ of the fusion algebra $\cA$ is canonically determined: these are exactly the non-zero elements of $\cA$ that cannot be expressed as the sum of two non-zero elements. The elements of $\cG$ are called the \emph{irreducible elements} of the fusion algebra $\cA$ and we sometimes write $\cG = \Irred \cA$. Notice that conjugation preserves irreducibility. The \emph{intrinsic group} $\grp(\cA)$ of the fusion algebra $\cA$ consists of the irreducible elements $x \in \cA$ such that $x \ox = e$. Equivalently, $x \in \cA$ belongs to the intrinsic group if and only if $x \ox$ is irreducible. It is easy to check that the intrinsic group of a fusion algebra is indeed a group. If $x,y \in \cA$, we sometimes say that \emph{$x$ is included in $y$}, if there exists a $z$ such that $y = x+z$.

A \emph{dimension function} on a fusion algebra $\cA$ is an additive, multiplicative, unital map $d : \cA \recht \R^+$ satisfying $d(\ox) = d(x)$ for all $x \in \cA$. Suppose $d$ is a dimension function on the fusion algebra $\cA$. Whenever $x \in \cA$ is non-zero, $e$ is included in $x \ox$ and so $d(x) \geq 1$. It then follows that $x \in \cA$ belongs to the intrinsic group of $\cA$ if and only if $d(x) = 1$. Moreover, if $x \in \cA$ is non-zero and not in the intrinsic group, the same reasoning yields $d(x) \geq \sqrt{2}$.

Two examples of fusion algebras with a dimension function arise as follows.
\begin{itemize}
\item Let $\Gamma$ be a group and define $\cA = \N[\Gamma]$. Define $d$ such that $d(s) = 1$ for all $s \in \Gamma$.
\item Let $G$ be a compact group and define the fusion algebra $\Rep(G)$ as the set of equivalence classes of finite dimensional unitary representations of $G$. The operations on $\Rep(G)$ are of course given by direct sum and tensor product of representations, while the dimension function $d$ is given by the ordinary Hilbert space dimension of the representation space.
\end{itemize}

We are interested in the fusion algebra $\FAlg(M)$ of a II$_1$ factor $M$. First of all, an $M$-$M$-bimodule $_M H_M$ is said to be of \emph{finite Jones index} if $\dim( _M H) < \infty$ and $\dim(H_M) < \infty$. We define $\FAlg(M)$ as the set of finite index $M$-$M$-bimodules modulo unitary equivalence.

Whenever $\psi : M \recht pM^n p$ is a finite index inclusion in the sense of Jones \cite{bc}, for some non-zero projection $p \in M^n$, we define the $M$-$M$-bimodule $H(\psi)$ on the Hilbert space $(\M_{1,n}(\C) \ot L^2(M))p$ with left and right module actions given by
$$a \cdot \xi := a \xi \quad\text{and}\quad \xi \cdot a = \xi \psi(a) \; .$$
Every finite index $M$-$M$-bimodule is unitarily equivalent with some $H(\psi)$. Moreover, given finite index inclusions $\psi : M \recht p M^n p$ and $\eta : M \recht q M^m q$, we have $H(\psi) \cong H(\eta)$ if and only if there exists a unitary $u \in p (\M_{n,m}(\C) \ot M)q$ satisfying $\psi(a) = u \eta(a) u^*$ for all $a \in M$.

Addition in $\FAlg(M)$ is given by the obvious direct sum of bimodules, while multiplication in $\FAlg(M)$ is given by the \emph{Connes tensor product of $M$-$M$-bimodules} that we recall now (see also V.Appendix~B in \cite{connes NCG}). Let $H$ be an $M$-$M$-bimodule. Define $\cH$ as the dense subspace of $H$ consisting of bounded vectors:
\[ \cH := \{\xi \in H \mid \; \exists c>0,\ \forall a \in M,\ \|\xi a\|_2 \leq c \|a\|_2\} \; .\] For all $\xi \in
\cH$ and $a \in M$ we define $L_\xi(a) = \xi a$. By definition this map extends to a bounded operator
$L_\xi: L^2(M) \recht H$. We set:
\[\langle \xi , \eta \rangle_M := L_\xi^* L_\eta \in M,\ \forall \xi,\eta \in \cH \; .\] It is easy to check
that this formula defines an $M$-valued scalar product on $\cH$. Then the Connes tensor product of the
$M$-$M$-bimodules $_MH_M$ and $_MK_M$ is defined as the separation and completion of the algebraic tensor product
$\cH \otimes_{\rm alg} K$ for the scalar product \[\langle a \ot \xi , b \ot \eta \rangle:= \langle \xi , \langle
a, b\rangle_M \eta \rangle \; . \] The Connes tensor product is denoted by $H \otimes_M K$ and it is an
$M$-$M$-bimodule: \[a \cdot (b \ot \xi) = ab \ot \xi\ \text{and}\ (b \ot \xi) \cdot a = b \ot(\xi
a) \; .\]
Note that $H(\psi) \ot_M H(\eta) = H((\id \ot \psi)\eta)$.

If $_M H_M \in \FAlg(M)$, the \emph{conjugate bimodule} $_M \overline{H}_M$ lives on the conjugate Hilbert space $\overline{H} = H^*$ with bimodule action given by
$$a \cdot \overline{\xi} = \overline{\xi a^*} \quad\text{and}\quad \overline{\xi} \cdot a = \overline{a^* \xi} \; .$$

The elements of the \emph{intrinsic group} $\grp(M)$ of $\FAlg(M)$ are exactly the bimodules $H(\pi)$, where $\pi : M \recht p M^n p$ is a $^*$-isomorphism. Denote by $\cF(M)$ the \emph{fundamental group} of $M$. We then get a short exact sequence $e \recht \Out(M) \recht \grp(M) \recht \cF(M) \recht e$, mapping $\si \in \Aut(M)$ to $H(\si) \in \grp(M)$ and mapping $H(\pi) \in \grp(M)$ to $\Tr(p)$.

The fusion algebra $\FAlg(M)$ has a natural dimension function: the dimension of $H(\psi)$ is defined as the square root of the \emph{minimal index} of $\psi(M) \subset pM^n p$. Since for an irreducible subfactor the minimal index equals the usual Jones index, the dimension function $d$ is given by
\[d( _M H_M) = \sqrt{\dim
(H_M)\; \dim ( _M H)} \; ,\] whenever $_M H_M$ is an \emph{irreducible} $M$-$M$-bimodule. We refer to \cite{PimPopa,Havet} for details.

\subsection{Minimal actions of compact groups and fusion algebras} \label{subsec.minimal}

A continuous action $G \actson M$ of a compact group $G$ on the II$_1$ factor $M$ is said to be \emph{minimal} if the map $G \recht \Aut(M)$ is injective and if $M \cap (M^G)' = \C 1$. Here, $M^G$ denotes the von Neumann algebra of $G$-fixed points in $M$.

Given such a minimal action $G \actson M$, set $N := M^G$. We get a canonical, dimension preserving, embedding $\Rep(G) \hookrightarrow \FAlg(N)$ of fusion algebras, defined as follows. Let $\pi : G \recht \cU(n)$ be an irreducible unitary representation of $G$. We choose a unitary $V_\pi \in \M_n(\C) \ot M$ satisfying
$$(\id \ot \si_g)(V_\pi) = V_\pi (\pi(g) \ot 1)$$
for all $g \in G$. We then define the finite index inclusion $$\psi_\pi : N \recht \M_n(\C) \ot N : \psi_\pi(a) = V_\pi (1 \ot a)V_\pi^* \; .$$
It is easily checked that the $N$-$N$-bimodule $H(\psi_\pi)$ is irreducible and, up to unitary equivalence, independent of the choice of $V_\pi$. The map $\pi \mapsto H(\psi_\pi)$ extends to an embedding $\Rep(G) \hookrightarrow \FAlg(N)$.

Also note that the coefficients of $V_\pi$ quasi-normalize $N$ and so, the inclusion $N \subset M$ is \emph{quasi-regular}.

\subsection{Freeness and free products of fusion algebras}

\begin{defi} \label{def.freeness}
Let $\cA$ be a fusion algebra and $\cA_i \subset \cA$ fusion subalgebras for $i = 1,2$. We say that $\cA_1$ and $\cA_2$ are \emph{free inside $\cA$} if every alternating product of irreducibles in $\cA_i \setminus \{e\}$, remains irreducible and different from $\{e\}$.
\end{defi}

Given fusion algebras $\cA_1$ and $\cA_2$, there is up to isomorphism a unique fusion algebra $\cA$ generated by copies of $\cA_1$ and $\cA_2$ that are free. We call this unique $\cA$ the \emph{free product} of $\cA_1$ and $\cA_2$ and denote it by $\cA_1 * \cA_2$. Of course, the free product can be constructed in a concrete way as follows: given $\cA_1$ and $\cA_2$, set $\cG_i = \Irred(\cA_i)$. Define $\cG$ as the set of words with letters alternatingly from $\cG_1 \setminus \{e\}$ and $\cG_2 \setminus \{e\}$. Denote the empty word as $e$. Then, $\cA_1 * \cA_2 = \N[\cG]$. The product on $\N[\cG]$ is the unique associative and distributive operation satisfying the following two conditions:
\begin{itemize}
\item The embeddings $\cA_i \hookrightarrow \N[\cG]$ are multiplicative.
\item If the last letter of the alternating word $x \in \cG$ and the first letter of the alternating word $y \in \cG$ belong to different fusion algebras $\cA_i$, the product of $x$ and $y$ is again irreducible and given by concatenation of $x$ and $y$.
\end{itemize}

Denote by $R$ the hyperfinite II$_1$ factor. It is a crucial ingredient of this paper that $\FAlg(R)$ is huge, in the sense that $\FAlg(R)$ contains many free fusion subalgebras. More precisely, it was shown in Theorem 5.1 of \cite{vaes} that countable fusion subalgebras of $\FAlg(R)$ can be made free by conjugating one of them with an automorphism of $R$ (see Theorem \ref{Gdelta} below). Note that the same result has first been proven for countable subgroups of $\Out(R)$ in \cite{IPP}. In both cases, the key ingredients come from \cite{P10}.

Let $M$ be a II$_1$ factor and $_M K_M \in \FAlg(M)$. Whenever $\al \in \Aut(M)$, we define the conjugation of $K$ by $\al$ as the bimodule $K^\al := H(\al^{-1}) \ot_M K \ot_M H(\al)$. Of course, $K^\al$ has $K$ as its underlying Hilbert space with new left and right module action given by $\xi \cdot_{\text{\rm new}} a = \xi \cdot_{\text{\rm old}} \al(a)$ and $a \cdot_{\text{\rm new}} \xi = \al(a) \cdot_{\text{\rm old}} \xi$.

\begin{theo}[Thm.\ 5.1 in \cite{vaes}] \label{Gdelta}
Let $R$ be the hyperfinite II$_1$ factor and $\cA_0,\cA_1$ two countable fusion
subalgebras of $\FAlg(R)$. Then, \[\{\alpha \in \Aut(R) \mid  \cA_0^{\alpha}\ \text{and}\
\cA_1\ \text{are free}\}\] is a $G_{\delta}$-dense subset of $\Aut(R)$.
\end{theo}

\subsection{Property (T) and relative property (T) for II$_1$ factors}

Property (T) for finite von Neumann algebras was defined by Connes and Jones in
\cite{T}: a II$_1$ factor $(N,\tau)$ has property (T) if and only if there exists $\epsilon >0$ and
a finite subset $F\subset N$ such that every $N$-$N$-bimodule $H$ that has a unit vector $\xi$ satisfying
$\|x \xi - \xi x\| \leq \epsilon,\ \forall x \in F$, actually has a non-zero $N$-central vector $\xi_0$, meaning that $x \xi_0 = \xi_0
x,\ \forall x \in N$.

Note that an ICC group $\Gamma$ has property (T) if and only if the II$_1$ factor $\cL(\Gamma)$ has property (T) in the sense of Connes and Jones.

Relative property (T) for inclusions $B \subset (N,\tau)$ of finite von Neumann algebras was defined by Popa in \cite{popa1}. We do not really use relative property (T) in this paper, just the trivial observation that $B \subset N$ has the relative property (T) if $B$ has itself property (T).

\section{Statement of the main result}

The main result that we prove is that any compact group $G$ can be realized as the
outer automorphism group of a type II$_1$ factor, see Theorem \ref{thm.exists}. A more precise theorem can be formulated as follows.

Note that any character $\om \in \Char(\Gamma)$ defines an automorphism $\theta_\om$ of any crossed product $N \rtimes \Gamma$ acting identically on $N$ and multiplying by $\om$ on $\Gamma$.

\begin{theo}\label{thm.thm}
Let $M_1$ be the hyperfinite II$_1$ factor and $G$ a compact group acting on $M_1$. Denote $N = M_1^G$, the von Neumann algebra of $G$-fixed points in $M_1$. Let $\Gamma$ be an ICC group acting on $N$. Denote $M_0 := N \rtimes \Gamma$. Assume that
\begin{enumerate}
\item the action $G \actson M_1$ is minimal,
\item the action $\Gamma \actson N$ is outer and $M_0$ has the property (T),
\item the natural images of $\Rep G \hookrightarrow \FAlg(N)$ and $\Aut(N \subset M_0) \overset{\text{\rm restr}}{\longrightarrow} \Out(N) \subset \FAlg(N)$ inside the fusion algebra $\FAlg(N)$, are free in the sense of Definition \ref{def.freeness}.
\end{enumerate}
Then, the homomorphism
$$\Char(\Gamma) \times G \recht \Aut(M_0 \underset{N}{*} {M_1}) : (\om,g) \mapsto \theta_\om * \si_g$$
induces an isomorphism $\Char(\Gamma) \times G \cong \Out(M_0 \underset{N}{*} {M_1})$.
\end{theo}

Combining Theorems \ref{Gdelta} and \ref{thm.thm}, we shall prove the following.

\begin{corol} \label{cor.cor}
Let $G$ be a compact, second countable group and $G \actson R$ a minimal action on the hyperfinite II$_1$ factor $R$. Then there exists an outer action of $\SL(3,\Z)$ on the fixed point algebra $R^G$, such that for $M$ given as the amalgamated free product $M = (R^G \rtimes \Gamma) \underset{R^G}{*} R$, the natural homomorphism
$$G \recht \Aut(M) : g \mapsto \id * \si_g$$
induces an isomorphism $G \cong \Out(M)$.
\end{corol}

Of course, every compact, second countable group $G$ admits a minimal action on the hyperfinite II$_1$ factor. A possible construction goes as follows: take an amenable ICC group $\Lambda$ and define the Bernoulli action crossed product $R = L^\infty\bigl(\prod_\Lambda (G,\text{\rm Haar})\bigr) \rtimes \Lambda$ with $G \actson R$ acting by diagonal left translation on $\prod_\Lambda (G,\text{\rm Haar})$ and trivially on $\Lambda$.

So, we immediately get the following result.

\begin{theo} \label{thm.exists}
Let $G$ be a compact, second countable group. There exists a type II$_1$ factor $M$ with $\Out(M) \cong G$.
\end{theo}

\section{A Kurosh automorphism theorem for fusion algebras}

An important ingredient in the analysis of all automorphisms of an amalgamated free product $M = M_0 *_N M_1$ as in Theorem \ref{thm.thm} above, is a generalization of the Kurosh automorphism theorem to automorphisms of free products of fusion algebras. We do not prove a general result, but a rather easy theorem sufficient for our purposes.

Recall that a group is said to be freely indecomposable if it cannot be written as a non-trivial free product.

\begin{theo} \label{Kurosh}
Let $\Gamma$ be a countable group non isomorphic to $\Z$ and freely indecomposable. Let
$\cA$ be an abelian fusion algebra with a dimension function $d$ and non-isomorphic to the group $\Z$.

Every dimension preserving automorphism $\al$ of $\N[\Gamma] * \cA$ is of the form $(\Ad u) \circ (\al_0 * \al_1)$ for some $u \in \Gamma * \grp(\cA)$, $\al_0 \in \Aut(\Gamma)$ and $\al_1$ a dimension preserving automorphism of $\cA$.
\end{theo}

\begin{proof}
Let $\al$ be a dimension preserving automorphism of $\N[\Gamma] * \cA$.

We denote $\Lambda = \grp(\cA)$, the intrinsic group of $\cA$, and $\Delta = \Gamma * \Lambda$, which is as well the intrinsic group of $\N[\Gamma] * \cA$. We also write $\cG = \Irred \cA$, which means that $\cA = \N[\cG]$. We may of course assume that $\cG \neq \Lambda$, because the group case of our theorem is covered by the classical Kurosh theorem. Finally, we set $\cGr = \cG \setminus \{e\}$ and $\Gammar = \Gamma \setminus \{e\}$. If $u \in \Delta$, we write $u^{-1}$ instead of $\overline{u}$.

{\bf Claim.} There exists $x \in \cG \setminus \Lambda$ and $u \in \Delta$ such that $\al(x) \in u (\cG \setminus \Lambda) u^{-1}$.

{\bf Proof of the claim.} Define $\lambda = \inf \{ d(x) \mid x \in \cG \setminus \Lambda \} \geq \sqrt{2}$. Take $x \in \cG \setminus \Lambda$ with $d(x) < \sqrt{2} \lambda$. Write $\al(x)$ as an alternating word in $\cGr$ and $\Gammar$. Suppose that in this expression of $\al(x)$, there appears twice a letter from $\cG \setminus \Lambda$. Then the dimension of these two letters is greater or equal than $\lambda \geq \sqrt{2}$, making $d(x) = d(\al(x)) \geq \sqrt{2} \lambda$~; a contradiction. So, we have shown that $\al(x) = u y v^{-1}$ with $y \in \cG \setminus \Lambda$ and $u,v \in \Delta$. We may assume that $u,v$ are either equal to $e$, either end with a letter from $\Gammar$. Expressing the commutation of $\al(x)$ and $\al(\ox)$, we find that $u y \overline{y} u^{-1} = v \overline{y} y v^{-1}$. Since $y \not\in \Lambda$, we find that $y \overline{y} \neq e$ and so $u = v$, proving the claim.

{\bf Observation 1.} If $x \in \Delta$, $y \in \cGr$ and $xy = yx$, then $x \in \Lambda$. This follows by analyzing reduced words in $\Gammar$ and $\cGr$.

Because of the claim and replacing $\al$ by $(\Ad u^{-1}) \circ \al$, we may from now on assume the existence of $x,y \in \cG \setminus \Lambda$ with $\al(x) = y$. Whenever $a \in \Lambda$, $\al(a)$ belongs to $\Delta$ and commutes with $y$. Observation~1 above implies that $\al(\Lambda) \subset \Lambda$. Similarly, $\al^{-1}(\Lambda) \subset \Lambda$ so that $\al(\Lambda) = \Lambda$. It follows that the restriction of $\al$ to $\Delta$ defines an automorphism of $\Gamma * \Lambda$ that globally preserves $\Lambda$. The classical Kurosh theorem implies that $\al(\Gamma) = \Gamma$.

{\bf Observation 2.} If $z \in \cGr$ and $\al(z) \in \Delta \cGr \Delta$, then actually $\al(z) \in \cG^\circ$. Indeed, write $\al(z) = u r v^{-1}$ for $r \in \cGr$ and $u,v \in \Delta$ either equal to $e$ or with their last letter in $\Gammar$. Writing out that $\al(z) = u r v^{-1}$ and $y = \al(x)$ commute, it follows that $u=v=e$. A similar observation holds for $\al^{-1}$.

It remains to prove that $\al(\cG) = \cG$. Assume the contrary and define
$$\delta = \inf \{ d(z) \mid z \in \cG, \quad (\; \al(z) \not\in \cG \;\;\text{or}\;\; \al^{-1}(z) \not\in \cG \; ) \; \} \; .$$
Take $z \in \cGr$ with $d(z) < \sqrt{2} \delta$ such that $\al(z) \not\in \cGr$ or $\al^{-1}(z) \not\in \cGr$. Assume that we are in the case $\al(z) \not\in \cGr$. By construction $\al(r), \al^{-1}(r) \in \cGr$ for every $r \in \cGr$ with $d(r) < \delta$. Write $\al(z)$ as an alternating word in $\cGr$ and $\Gammar$. By observation~2, the expression for $\al(z)$ contains at least twice a letter from $\cG \setminus \Lambda$. Hence every letter in the expression for $\al(z)$ has dimension strictly smaller than $\delta$. Applying $\al^{-1}$ and using the fact that $\al^{-1}(\Gammar) = \Gammar$, we have written $z \in \cGr$ as an alternating word in $\Gammar$ and $\cGr$ with more than $2$ letters; a contradiction.
\end{proof}

\section{Proof of the main results}

Before proving Theorem \ref{thm.thm} and Corollary \ref{cor.cor}, we state the following lemma. It is a consequence of Lemma 8.4 in \cite{IPP} (see also Props.\ 3.3 and 3.5 in \cite{vaes}).

\begin{lemm}\label{IPP2}
Let $\Gamma_0,\ \Gamma_1$ be ICC groups acting outerly on the II$_1$ factors $A_0$ and $A_1$
respectively. Set $M := A_0 \rtimes \Gamma_0$ and suppose that $\alpha: A_0 \rtimes \Gamma_0 \recht A_1
\rtimes \Gamma_1$ is an isomorphism such that $\alpha(A_0) \prec_M A_1$ and $A_1 \prec_M \alpha(A_0)$.

Then, there exists a unitary $u \in \cU(M)$ such that $u \alpha(A_0)u^* = A_1$.
\end{lemm}

A first step in the proof of Theorem \ref{thm.thm} is the following lemma. The crucial ingredients of its proof are Theorems 1.2.1 and 4.3 in \cite{IPP} (see also Thms.\ 4.6 and 5.6 in \cite{houdayer} for alternative proofs).

\begin{lemm}\label{lem.firststep}
Suppose that the assumptions of Theorem \ref{thm.thm} are fulfilled. Set $M = M_0 *_N M_1$. For every $\al \in \Aut(M)$, there exists $u \in \cU(M)$ such that $(\Ad u) \circ \al \in \Aut(N \subset M)$.
\end{lemm}

Note that in fact assumption~3 in Theorem \ref{thm.thm} will not be used in the proof of this lemma.

\begin{proof}
By 4.3 in \cite{IPP} and because $M_0$ has property (T), there exists $i \in \{0,1\}$ such that $\al(M_0) \prec_M M_i$. Since $M_1$ is hyperfinite, it follows that $i = 0$. So, we can take a projection $p \in M_0^n$, a non-zero partial isometry $v \in \M_{1,n}(\C) \ot M$ and a unital $^*$-homomorphism $\psi : M_0 \recht p M_0^n p$ satisfying $\al(a) v = v \psi(a)$ for all $a \in M_0$. Since $\psi(M_0)$ has property (T), we know that $\psi(M_0) \not\prec_{M_0^n} N^n$. By 1.2.1 in \cite{IPP}, it follows that $\psi(M_0)' \cap p M^n p \subset p M_0^n p$. In particular, $v^* v \in p M_0^n p$. So, we may assume that $p = v^* v$. Since also $\al(M_0)' \cap M = M_0' \cap M = \C 1$, we have $v v^* = 1$. Factoriality of $M_0$ now allows to assume that $v \in \cU(M)$ and $v^* \al(M_0) v \subset M_0$.

Applying the same reasoning to $\al^{-1}$, we also get a unitary $w \in \cU(M)$ satisfying $w^* M_0 w \subset \al(M_0)$. It follows that $(wv)^* M_0 (wv) \subset M_0$. Another application of 1.2.1 in \cite{IPP} implies that $wv \in M_0$. But then all the inclusions $(wv)^* M_0 (wv) \subset v^* \al(M_0) v \subset M_0$ are equalities. So, after a unitary conjugacy, we may assume that $\al(M_0) = M_0$.

Quasi-regularity of $N \subset M$ combined with 1.2.1 in \cite{IPP}, implies that $\al(N) \prec_{M_0} N$. Similarly $N \prec_{M_0} \al(N)$. The lemma then follows by applying Lemma \ref{IPP2}.
\end{proof}

\begin{proof}[\bf Proof of Theorem \ref{thm.thm}]
We still denote $M = M_0 *_N M_1$. Let $\al \in \Aut(M)$. We prove below that after a unitary conjugacy of $\al$, one has $\al(a) = a$ for all $a \in N$ and $\al(M_i) = M_i$ for $i \in \{0,1\}$. This then implies that $\al|_{M_0} = \theta_\om$ for some $\om \in \Char(\Gamma)$ and $\al|_{M_1} = \si_g$ for some $g \in G$. Hence, it implies the surjectivity of the homomorphism $\Char(\Gamma) \times G \recht \Out(M)$. The injectivity of this homomorphism follows from the irreducibility $N' \cap M = \C 1$ that we prove now as the consequence of more general needed considerations.

Let $\cI$ be a complete set of inequivalent irreducible unitary representations of $G$. For every $\pi \in \cI$, choose a unitary $V_\pi \in \B(H_\pi) \ot M_1$ satisfying $(\id \ot \si_g)(V_\pi) = V_\pi (\pi(g) \ot 1)$. Define $K_0(\pi) \subset M_1$ as the linear span of
$$(\xi^* \ot a)V_\pi (\eta \ot 1) \;\; , \;\; \xi,\eta \in H_\pi \;\; , \;\; a \in N \; .$$
It follows that the closure $K(\pi)$ of $K_0(\pi)$ is a finite index $N$-$N$-subbimodule of $L^2(M_1)$. Denote by
$$\Psi : \Rep G \hookrightarrow \FAlg(N)$$
the embedding defined in Subsection \ref{subsec.minimal}. It follows that $K(\pi) \cong (\dim \pi) \cdot \Psi(\pi)$. Moreover, we have the following orthogonal decomposition of $L^2(M_1)$.
$$L^2(M_1) = \bigoplus_{\pi \in \cI} K(\pi) \; .$$
In the same way, we define for every $s \in \Gamma$, the subspace $H_0(s) \subset M_0$ given by $H_0(s) = N u_s$, with closure $H(s) \subset L^2(M_0)$.

Whenever $w= s_0 \pi_1 s_1 \cdots s_{n-1} \pi_n s_n$ is an alternating word in $\Gamma \setminus \{e\}$ and $\cI \setminus \{e\}$, we define the $N$-$N$-subbimodule $H(w) \subset L^2(M)$ as the closure of $H_0(s_0) K_0(\pi_1) \cdots K_0(\pi_n) H_0(s_n)$. One then obtains the orthogonal decomposition of $L^2(M)$ given by
\begin{equation} \label{eq.crucialone}
L^2(M) = \bigoplus_{w \; \text{\rm alternating word}} H(w) \; .
\end{equation}
Moreover, as $N$-$N$-bimodules, we get the unitary equivalences
\begin{equation} \label{eq.crucialtwo}
\begin{split}
H(w) & \cong H(s_0) \ot_N K(\pi_1) \ot_N \cdots \ot_N K(\pi_n) \ot_N H(s_n) \\ & \cong \;\text{a multiple of}\; H(s_0) \ot_N \Psi(\pi_1) \ot_N \cdots \ot_N \Psi(\pi_n) \ot_N H(s_n) \; .
\end{split}
\end{equation}
Denote by $\FAlg(N \subset M)$ the fusion subalgebra of $\FAlg(N)$ generated by the finite index $N$-$N$-subbimodules of $L^2(M)$. Using assumption~3 in Theorem \ref{thm.thm} and \eqref{eq.crucialone} and \eqref{eq.crucialtwo} above, $\Psi$ extends to an isomorphism $\Psi : \N[\Gamma] * \Rep(G) \recht \FAlg(N \subset M)$ such that $H(w)$ is a multiple of $\Psi(w)$ for every alternating word $w$. We see in particular that $L^2(M)$ contains only once the trivial $N$-$N$-bimodule (for $w$ the empty alternating word and $H(w) = L^2(N)$). This means that $N' \cap M = \C 1$ as was needed above.

Let $\Char G \subset \cI$ be the subset of $\cI$ consisting of one-dimensional unitary representations of $G$. Then, $\Char G$ is as well the intrinsic group of the fusion algebra $\Rep G$. Whenever $\pi \in \Char G$, we have $V_\pi \in \cU(M_1)$ and $V_\pi$ normalizes $N$. Whenever $w \in \Gamma * \Char G$, write $w = s_0 \pi_1 s_1 \cdots s_{n-1} \pi_n s_n$ as an alternating word in $\Gamma \setminus \{e\}$ and $\Char G \setminus \{e\}$ and define the unitary $U(w):= u_{s_0} V_{\pi_1} \cdots V_{\pi_n} u_{s_n}$ normalizing $N$.

We are now ready to complete the proof of the theorem. So, let $\al \in \Aut(M)$. We have to prove that after a unitary conjugacy of $\al$, one has $\al(a) = a$ for all $a \in N$ and $\al(M_i) = M_i$ for $i \in \{0,1\}$. By Lemma \ref{lem.firststep}, we may assume that $\al(N) = N$. But then, the conjugation map $K \mapsto K^\al$ defines an automorphism of the fusion subalgebra $\FAlg(N \subset M)$ of $\FAlg(N)$. Define the automorphism $\eta$ of $\N[\Gamma] * \Rep(G)$ such that $\Psi(\eta(w)) = (\Psi(w))^\al$ in $\FAlg(N)$ for all $w \in \N[\Gamma] * \Rep(G)$. By Theorem \ref{Kurosh}, we find an element $v$ in $\Gamma * \Char(G)$ such that $(\Ad v) \circ \eta$ globally preserves $\Gamma$ and $\Rep(G)$.

Replacing $\al$ by $(\Ad U(v)) \circ \al$, we may assume that $\eta$ preserves globally $\Gamma$ and $\Rep(G)$. The equality $\al(N) = N$ remains true. The restrictions of $\eta$ yield an automorphism of the group $\Gamma$ and a permutation of $\cI$ respecting the fusion rules. Moreover, we have $K(\pi)^\al \cong K(\eta(\pi))$ for every $\pi \in \cI$ and $H(s)^\al \cong H(\eta(s))$ for every $s \in \Gamma$. Choose $s \in \Gamma$. Note that $H(s)^\al$ is isomorphic as an $N$-$N$-bimodule with the closure of $\al^{-1}(N u_s)$ in $_N L^2(M)_N$. Since the $N$-$N$-bimodule $H(\eta(s))$ appears with multiplicity $1$ in the decomposition \eqref{eq.crucialone} of $L^2(M)$, we conclude that the closure of $\al^{-1}(N u_s)$ inside $L^2(M)$ equals $H(\eta(s))$ for all $s \in \Gamma$. It follows that $\al(M_0) = M_0$. A similar reasoning shows that $\al(M_1) = M_1$.

Since $\al \in \Aut(N \subset M_0)$, assumption~3 in Theorem \ref{thm.thm} implies that the $N$-$N$-bimodule $H(\al|_N)$ is free with respect to $\Psi(\Rep G)$ inside $\FAlg(N)$. But the formula $K(\pi)^\al \cong K(\eta(\pi))$ means that $H(\al|_N)$ normalizes $\Psi(\Rep G)$. Both statements can only be true at the same time if $H(\al|_N)$ is the trivial $N$-$N$-module. So, $\al|_N$ is an inner automorphism of $N$ and we are done.
\end{proof}

\begin{proof}[\bf Proof of Corollary \ref{cor.cor}]
It suffices to give an example of an outer action of $\Gamma=\SL(3,\Z)$ on the hyperfinite II$_1$ factor $N$ such that $N \rtimes \Gamma$ has property (T). Indeed, starting from a minimal action of $G$ on the hyperfinite II$_1$ factor $M_1=R$, set $N = M_1^G$ and take an outer action of $\Gamma$ on $N$ such that $M_0:=N \rtimes \Gamma$ has property (T). By Theorem 4.4 in \cite{popa1}, $\Aut(N \subset N\rtimes \Gamma)/ \Ad \cU(N)$ is a countable group. By Theorem \ref{Gdelta}, we can take an automorphism $\al \in \Aut(R)$ and replace $\Gamma \subset \Aut(N)$ by $\al \Gamma \al^{-1}$ in such a way that all conditions of Theorem \ref{thm.thm} are fulfilled. Since $\Char \Gamma = \{e\}$, the corollary then follows from Theorem \ref{thm.thm}.

Take $\Gamma_1: = \SL(3,\Z) \ltimes (\Z^3 \oplus \Z^3)$ with the action of $\SL(3,\Z)$ on $\Z^3 \oplus \Z^3$ given
$$ A \cdot(x,y):= (Ax, (A^{-1})^ty) \; .$$ Note that $\Gamma_1$ is a property (T) group.
Take $k \in \R \setminus 2 \pi \Q$ and define the non degenerate 2-cocycle $\Omega \in Z^2(\Z^3 \oplus \Z^3,
S^1)$ by the formula
$$\Omega \bigl( (x,y) ; (x',y') \bigr):=e^{\mathbf i k (\langle x, y' \rangle -
\langle y,x'\rangle)}$$
where $\langle \cdot , \cdot \rangle$ is the standard scalar product on $\Z^3$. The
$2$-cocycle $\Omega$ is $\SL(3,\Z)$-invariant and hence extends to a 2-cocycle $\tilde{\Omega} \in
Z^2(\Gamma_1, S^1)$ by the formula
$$\tilde{\Omega} \bigl( ( (x,y), A ); ( (x',y'), B ) \bigr): =
\Omega \bigl( (x,y) ; A \cdot (x', y') \bigr) \; .$$
The twisted group von Neumann algebra $\cL_{\tilde{\Omega}}(\Gamma_1)$ still has property (T) and can be regarded as well as $\cL_\Omega(\Z^3 \oplus \Z^3) \rtimes \SL(3,\Z)$. Since $\cL_\Omega(\Z^3 \oplus \Z^3)$ is the hyperfinite II$_1$ factor, we are done.
\end{proof}

\end{document}